\documentclass{amsart}

\usepackage{amssymb}
\usepackage[all]{xy}
\usepackage{hyperref}

\newtheorem{thm}{Theorem}

\newtheorem{lem}[thm]{Lemma}

\newtheorem{prop}[thm]{Proposition}

\theoremstyle{definition}

\newtheorem{say}[thm]{}
\newtheorem{exmp}[thm]{Example}

\newtheorem*{ack}{Acknowledgments}      

\newtheorem{defn-thm}[thm]{Definition--Theorem}  
\newtheorem{defn-lem}[thm]{Definition--Lemma}  

\theoremstyle{remark}
\newtheorem{claim}[thm]{Claim}

\renewcommand{\c}[0]{{\mathbb C}}

\renewcommand{\r}[0]{{\mathbb R}}

\newcommand{\s}[0]{{\mathbb S}}

\newcommand{\qtq}[1]{\quad\mbox{#1}\quad}

\newcommand{\tsum}[0]{\textstyle{\sum}}

\def\loccoh#1.#2.#3.#4.{H^{#1}_{#2}(#3,#4)}

\DeclareMathAlphabet{\mathchanc}{OT1}{pzc}%
                                {m}{it}

\usepackage[all]{xy}\xyoption{dvips}

\begin{document}
\bibliographystyle{amsalpha}

 \title{Checking real analyticity on surfaces}
 \author{Jacek Bochnak, J\'anos Koll\'ar \and Wojciech Kucharz}

\address{Jacek Bochnak, Le Pont de l'Etang 8 1323
Romainm\^otier, Switzerland}
\email{jack3137@gmail.com}

\address{J\'anos Koll\'ar,  Princeton University, Princeton NJ 08544-1000 USA}
\email{kollar@math.princeton.edu}

\address{Wojciech Kucharz, Institute of Mathematics, Faculty of Mathematics and Computer
Science, Jagiellonian University, \L{}ojasiewicza 6, 30-348
Krak\'ow, Poland}
\email{Wojciech.Kucharz@im.uj.edu.pl}

\begin{abstract} 
We prove that a real-valued function (that is not assumed to be  continuous) on a real analytic manifold is
analytic whenever all its restrictions to   analytic
submanifolds homeomorphic to $\s^2$ are analytic. 
This is a real analog  for the
classical theorem of Hartogs that a function on a complex manifold   is complex analytic iff it is 
complex  analytic when restricted to any complex curve.
\end{abstract}

 \maketitle

By a theorem of Hartogs, a function $f:\c^n\to \c$ (that is not assumed to be  continuous) is analytic iff it is 
analytic when restricted to any translate of the coordinate axes; see 
\cite{har} or  \cite[p.140]{bib4}.
An analogous claim does not hold for real analytic functions; 
see Example~\ref{2.examples}  for a simple non-continuous function
$f\colon \r^3\to \r$ that is real analytic  on every translate of the coordinate hyperplanes and a more complicated non-continuous function
$g\colon \r^3\to \r$ that is real analytic  on every smooth analytic curve in $\r^3$. 

Bochnak and Siciak   proved  in 1971  \cite{bib2, bib3} that 
a function $f:\r^n\to \r$  is real analytic  iff its  restriction to any 2-plane is  analytic. They also conjectured that a similar result holds on any 
real analytic manifold using  restrictions to $2$-dimensional compact analytic
submanifolds.

The aim  of this note is to prove  a stronger variant of this,   using an improvement of \cite[Thm.1]{bib3} that uses fewer 2-planes.

\begin{thm}\label{m.0.thm}  Let $M$ be a real analytic manifold of dimension $n\geq 3$ and
$f:M\to \r$ a (not necessarily continuous) function. 
Assume that  $f|_N$   is real analytic for every  real analytic submanifold $N\subset M$ that is homeomorphic to the 2-sphere $\s^2$.
 Then $f$ is 
 real analytic on $M$.
\end{thm}

Working on local coordinate charts, this is an immediate consequence of the following more precise version. 

\begin{thm}\label{m.thm}  Let $B\subset \r^n$ be the  open unit ball and
$f:B\to \r$ a (not necessarily continuous) function. 
Assume that  $n\geq 3$ and the restriction  $f|_{S^2}$  is real analytic for every 2-sphere $S^2\subset B$  passing through the origin. Then $f$ is 
 real analytic.
\end{thm}

Here we use 2-sphere  in the most restrictive sense, that is,
the intersection of an $(n-1)$-sphere
$\bigl (\tsum_i x_i^2=\tsum_i c_ix_i\bigr)$ with a  vector subspace of dimension 3.
\medskip

{\it Terminology.}
From now on  by a function we mean an arbitrary real-valued function and real analytic is shortened to analytic. To simplify notation,  we let  $g|_Y$ denote the restriction of $g$ to $Y\cap (\mbox{domain of $g$})$.

\medskip

Before we start the proof we need some preparation.
A subset  $C\subset \r^n$ is a {\it cone} if   $\r C= C$. (This is slightly more convenient for us than  the more usual variant  $\r^+ C= C$.) 
A cone $C$ is called open if $C\setminus \{0\}$ is open in $\r^n\setminus \{0\}$. A cone-neighborhood of a set $W$ is an open cone that contains $W$.

We will work both with vector subspaces  $V\subset \r^n$ and affine-linear
subspaces  $Q\subset \r^n$.  A property ${\mathcal P}$  holds for all {\it vector subspaces  near} $V$
if there is a    cone-neighborhood $C\supset V$ such that  ${\mathcal P}$  holds for all vector subspaces (of the same dimension) contained in $C$. 
A property ${\mathcal P}$  holds for all {\it affine-linear subspaces  near} $Q$ 
if there is a choice of the origin $0\in Q$, a  cone-neighborhood $C\supset Q$ and $\epsilon>0$ such that  ${\mathcal P}$  holds for all  affine-linear subspaces (of the same dimension) contained in $C\cup B(\epsilon)$, where $B(\epsilon)$ denotes the ball  of radius $\epsilon$ centered at $0$.
 For a cone $W$ and $\epsilon>0$ we write
$W(\epsilon):=W\cap B(\epsilon)$.
We need  the following obvious fact.

\begin{claim} \label{claim.1} 
Let $C\subset \r^n$ be an open cone and 
 $V\subset C$  a vector 2-plane and let $C_V\subset C$ be the union of all  vector 2-planes   $V'\subset C$ for which $\dim (V'\cap V)\geq 1$. 
Then $C_V$  is an open subcone of $C$ that contains $V$.\qed
\end{claim}
By Lagrange interpolation, a 1-variable polynomial of degree $d$ is 
 uniquely determined by   $d+1$ of its values. Equivalently, 
a 2-variable homogeneous polynomial of degree $d$ on $\r^2$ is uniquely determined by  its restriction to $d+1$ vector lines. 

More generally, the following is easy to prove by induction on $d$.

\begin{claim} \label{interp.say.1}
Let $H_0,\dots, H_d\subset \r^n$ be vector hyperplanes such that the intersection of any $r$ of them has codimension $r$ for every $r\leq n$.
For $i=0,\dots, d$ let $g_i$ be a degree $d$ homogeneous polynomial on $H_i$ 
such that $g_i|_{H_i\cap H_j}=g_j|_{H_i\cap H_j}$ for all $i,j$.
Then there is a unique degree $d$ homogeneous polynomial $g$ on  $\r^n$ such that $g|_{H_i}=g_i$ for every $i$.\qed
\end{claim}

\begin{lem}  \label{lem.1} Let  $V\subset \r^n$ be a vector 2-plane,  $C\supset V$  a  cone-neighborhood and $C_V\subset C$ as in  Claim~\ref{claim.1}.
Let $\sigma:C\to \r$ be a  function
such that $\sigma|_{V'}$ is a degree $d$ homogeneous polynomial 
for every vector 2-plane $V'\subset C_V$.
Then there is a unique degree $d$ homogeneous polynomial $g$ on $\r^n$ 
 such that 
$g|_{C_V}=\sigma|_{C_V}$. 
\end{lem}

Proof.  By induction on $n$, starting with $n=2$ when there is nothing to prove. For $n\geq 3$ pick 
 any  $V'\subset C_V$ such that  $\dim (V'\cap V)= 1$ and  let $H$ be  a vector  hyperplane that contains $V$ but does not contain $V'$. For every $i=0,\dots, d$ let $V_i\subset H_i$ be a general perturbation of  $V\subset H$  such that $V_i\subset C_V$,  the intersections
$V'\cap V_i$ are distinct lines and $V'\cap H_i=V'\cap V_i$.

By induction, there are degree $d$ homogeneous polynomials $g_i$ defined on
$H_i$ such that $g_i$ agrees with $\sigma$ on 
$H_i\cap C_V$.  By Claim~\ref{interp.say.1} there is thus a 
 unique degree $d$ homogeneous polynomial $g$ on  $\r^n$ such that $g|_{H_i}=g_i$ for every $i$.
 
 Thus
$g|_{V'}$ and $\sigma|_{V'}$ are 2 homogeneous polynomials of degree $d$
that agree on the $d+1$ lines  $V'\cap H_i$. So
$g|_{V'}=\sigma|_{V'}$. \qed

\begin{lem} \label{thm.1}
Let $V\subset \r^n$ be a vector 2-plane, $C\subset \r^n$ a cone-neighborhood of $V$   and $C_V\subset C$ as in  Claim~\ref{claim.1}.
Let  $f:C\to \r$ be a  function  such that,
for every vector  2-plane $V'\subset C_V$, there is an $\epsilon=\epsilon(V')>0$ such that 
the restriction of $f$ to $V'(\epsilon)$ is   analytic.
Then there is an analytic function $T$ in a neighborhood of $0\in \r^n$ 
 such that,
for every vector line  $L\subset C_V$, the restrictions
$f|_L$ and $T|_L$ agree  in a neigborhood of the origin in $L$.
\end{lem}

Proof. 
For any vector line $L\subset C_V$ the function $f|_{L}$ is  analytic at the origin. Thus, for $0\neq p\in L$  we can define
$$
\tau_r(p):=\tfrac{d^rf}{dt^r}(tp)\bigr|_{t=0}\qtq{and set}\tau_r(0)=0.
$$
By assumption $f|_{V'}$ is  analytic at the origin for every vector  2-plane $V'\subset C_V$, hence $\tau_r|_{V'}$ is a  degree $r$ homogeneous polynomial.
Furthermore, $f|_{L}=\tsum_r \tfrac1{r!}\tau_r$ near the origin.
By Lemma~\ref{lem.1}  there is a  degree $r$ homogeneous polynomial  $T_r$ on $\r^n$  that agrees with $\tau_r$ on    $C_V\subset C$. 

Consider the series of homogeneous polynomials
$$
T:=\tsum_r  \tfrac1{r!} T_r.
\eqno{(\ref{thm.1}.1)}
$$
As we noted,  this sum converges in some neighborhood of the origin
$0\in L(\epsilon_L)\subset L$ for every vector line $L\subset C_V$. 
Therefore  the series (\ref{thm.1}.1) defines an analytic function $T$ in a neighborhood  $B(\eta)\ni 0$ by  \cite[Lem.3]{bib3}.  Clearly $T$ has the required properties.
\qed

\medskip
{\it Remark.} It should be noted that the convergence of a series of homogeneous polynomials is quite different from the usual convergence of power series. In particular, if  the series (\ref{thm.1}.1)
 absolutely converges at a point  $(p_1,\dots, p_n)$, it does not imply that
it also converges at all points $(x_1,\dots, x_n)$ with  $|x_i|< |p_i|$. 
Thus \cite[Lem.3]{bib3} is  quite a subtle tool.
\medskip

Now we are ready to prove  the following stronger form of  \cite[Thm.1]{bib3}.

\begin{thm}  \label{m.thm.2}  Let $ U\subset \r^n$ be an open set and
$f:U\to \r$ a  function. 
Let $Q\subset \r^n$ be an affine 2-plane and assume that
$f|_{Q'\cap U}$ is  analytic for every affine  2-plane $Q'$ near   $Q$.
Then $f$ is  analytic in a neighborhood of $Q\cap U$. 
\end{thm}

Proof. The question is local, so we may assume that $0\in Q\cap U$ and work near it.  
Lemma~\ref{thm.1} gives a function $T$ that is 
 analytic in a ball  $B(\eta)$, and  a  cone-neighborhood $C$ of $Q$, such that  the restrictions
$f|_L$ and $T|_L$ agree  in a neigborhood of the origin
for every vector line  $L\subset C$. 
We may assume that  $B(\eta)\subset U$, thus $f$ and $T$ are both defined and  analytic on $L(\eta)$. Hence  $f$ and $T$ agree on $C(\eta)=C\cap B(\eta)$.

It remains to show that 
$T=f$ in a neighborhood of $0$. 
Fix an affine  line $\ell\subset Q$ not passing through the origin such that $\ell\cap  C(\eta)$ is not empty.
If $p\in \r^n$ is close to $0$ then the affine 2-plane $Q_p$ obtained as the span of $ p$ and $ \ell$ is close to $Q$. Thus   $T$ and $f$ restrict to   analytic functions  on the convex open set  $Q_p\cap B(\eta)$ that agree on the nonempty open subset
 $Q_p\cap  C(\eta)$. Hence they agree everywhere in a neigborhood of $Q\cap U$.  \qed

\begin{say}[Proof of Theorem~\ref{m.thm}]\label{m.thm.pf}
  Inversion
$$
\mu: (x_1,\dots, x_n)\mapsto \bigl(\tfrac{x_1}{\sum x_i^2},\dots, \tfrac{x_n}{\sum x_i^2}\bigr)
$$
maps the punctured unit ball $\bigl(0<\tsum x_i^2<1\bigr)$ to the outside of the unit ball
$\bigl(\tsum x_i^2>1\bigr)$ and  gives a one-to-one correspondence between the 2-spheres $S^2\subset B(1)$ that pass  through
the origin and   2-planes contained in the outside of the unit ball.  

For every $p\in \bigl(\tsum x_i^2>1\bigr)$ there is an affine 2-plane
$p\in Q_p\subset \bigl(\tsum x_i^2>1\bigr)$ and every affine 2-plane near $Q_p$  is also contained in $\bigl(\tsum x_i^2>1\bigr)$.
Thus $f\circ \mu^{-1}$ is  analytic by  Theorem~\ref{m.thm.2}, so $f$ is  analytic, except possibly at the origin. 

To see what happens at the origin, we use inversion 
$\mu_{\mathbf p}({\mathbf x}):=\mu({\mathbf x}-{\mathbf p})+{\mathbf p}$ centered at another point $0\neq {\mathbf p}\in B(1)$. This maps $B(1)\setminus\{{\mathbf p}\}$ to an open set  $U_{\mathbf p}\subset \r^n$. The 2-planes that pass through $\mu_{\mathbf p}(0)$ and are contained in  $U_{\mathbf p}$ are in
one-to-one correspondence with  the 2-spheres $S^2\subset B(1)$  passing  through both
the origin and ${\mathbf p}$. Thus $f\circ \mu^{-1}_{\mathbf p}$ is  analytic on these planes.
We already know that $f\circ \mu^{-1}_{\mathbf p}$ is  analytic away from $\mu_{\mathbf p}(0)$,
thus it is also  analytic at $\mu_{\mathbf p}(0)$  by  Theorem~\ref{m.thm.2}. 
So $f$ is also  analytic at the origin. \qed
\end{say}

The following is another variant of  Theorem~\ref{m.thm}.

\begin{prop} \label{cor.4}  
Let $M$ be an analytic manifold of dimension $n$ and 
$f:M\to \r$ a  function. 
Let $h:M\to \r$ be an analytic function that has an isolated 0 at a point $p\in M$ and let $x_1,\dots, x_n$  local coordinates at $p$. 
Assume that $f|_S$ is  analytic for every compact, 2-dimensional, smooth submanifold $S\subset M$ of the form  
$$
S:=\bigl(h-\ell_1=\ell_2=\cdots=\ell_{n-2}=0\bigr),
$$
where the $\ell_i:=\sum_j a_{ij}x_j$ are linear.

Then $f$ is 
  analytic in a punctured neighborhood of  $p\in M$. 
\end{prop}

Proof. We argue as in Paragraph~\ref{m.thm.pf}, but use the map
$ \bigl(\tfrac{x_1}{h},\dots, \tfrac{x_n}{h}\bigr)$,
which sends the submanifolds $S\subset M$ to affine 2-planes in $\r^n$.\qed

\begin{exmp}\label{2.examples}
The function defined by 
$$
f:=\frac{x_1\cdots x_n}{x_1^{2n}+\cdots+ x_n^{2n}}  \qtq{for}
(x_1,\dots, x_n)\neq (0,\dots,  0) \qtq{and}  f(0,\dots,  0)=0
$$
is  analytic  on every translate of the coordinate hyperplanes, but  not even bounded  at the origin. Thus we definitely need more planes than in the complex case.

Let $g \colon \r^3 \to \r$ be the function defined by
$$
g(x,y,z) = \frac{x^8 + y(x^2 - y^3)^2 +z^4}{x^{10} + (x^2 - y^3)^2 +z^2}
\qtq{for}  (x,y,z) \neq (0,0,0) \qtq{and}
g(0,0,0)=0.
$$
Then $f$ is analytic  on every nonsingular
analytic  curve in $\r^3$, but $f$ is not even
continuous at $(0,0,0)$.  See  also \cite{bmp} for an even stronger example.
\end{exmp}

\begin{ack}
JB and WK thank the Mathematisches Forschungsinstitut Oberwolfach for excellent
working conditions during their stay within the Research in Pairs
Programme. Partial  financial support  to JK   was provided  by  the NSF (USA) under grant number
 DMS-1362960. Partial financial support for WK was provided
by the National Science Center (Poland) under grant number
2014/15/B/ST1/00046.
\end{ack}


\begin{thebibliography}{BM48}



\bibitem[BMP91]{bmp} E.~Bierstone, P.~D.~Milman  and A.~Parusi{\'n}ski,
A function which is arc-analytic but not continuous,
Proc. Amer. Math. Soc.  113  (1991)  419--423.


\bibitem[BM48]{bib4} S.~Bochner and W.~Martin, Several Complex Variables,
Princeton University Press, Princeton (1948).



\bibitem[BS71]{bib2} J.~Bochnak and J.~Siciak,  Analytic
functions in topological vector spaces, mimeographed preprint, IHES  (1971).

\bibitem[BS18]{bib3} J.~Bochnak and J.~Siciak, A characterization of analytic
functions of several real variables, Ann. Polon. Math. (online 2018),
DOI: 10.4064/ap180119-26-3.


\bibitem[Har1906]{har}
Fritz~Hartogs, Zur {T}heorie der analytischen {F}unktionen mehrerer
              unabh\"{a}ngiger {V}er\"{a}nderlichen, insbesondere \"{u}ber 
die  {D}arstellung derselben durch {R}eihen, welche nach 
{P}otenzen   einer {V}er\"{a}nderlichen fortschreiten,
 Math. Ann.  62  (1906), 1--88.


\end{thebibliography}
\end{document}